\documentclass[pdflatex,sn-mathphys-ay]{sn-jnl}

\usepackage{amsmath, amssymb, bm, graphicx}

\usepackage{hyperref,doi}

\newcommand{\e}{\mathrm e}
\newcommand{\pfrac}[2]{\frac{\partial #1}{\partial #2}}

\begin{document}

\title{Wavespeed selection and interstitial gap formation in an acid-mediated cancer invasion model}

\author[1]{\fnm{Yuhui} \sur{Chen}}

\author*[1]{\fnm{Michael} \sur{Dallaston}}\email{michael.dallaston@qut.edu.au}

\affil[1]{\orgdiv{School of Mathematical Sciences}, \orgname{Queensland University of Technology}, \orgaddress{\street{2 George St}, \city{Brisbane}, \postcode{4000}, \state{QLD}, \country{Australia}}}

\abstract{
We consider a two-component reaction-diffusion system that has previously been developed to model invasion of cells into a resident cell population.  The system is an idealised version of models of tumour growth in which tumour cells degrade the surrounding tissue by increasing the acidity of the local environment.
By numerically computing families of travelling wave solutions to this problem, we observe that a general initial condition with either compact support, or sufficiently large exponential decay in the far field, tends to the travelling wave solution that has the largest possible decay at its front.  Initial conditions with sufficiently slow exponential decay tend to those travelling wave solutions that have the same exponential decay as their initial conditions.
We also show that in the limit that the (nondimensional) degradation rate of resident cells is large, the system has similar asymptotic structure as previously observed in perturbed Fisher--KPP models.  The asymptotic analysis in this limit explains the formation of an interstitial gap (a region between the invading and receding fronts, in which both cell populations are small), the width of which is logarithmically large in the limit of large degradation rate.
These results show that the general mechanism behind the formation of the interstitial gap in reaction-diffusion tumour models is connected to perturbations of the Fisher-KPP system.  Biologically, this implies that order of magnitude difference in degradation rate is required to produce appreciably different gap sizes, while the velocity of the invading front is largely determined by the Fisher-KPP velocity, and only very weakly affected by the presence of the interstitial gap.}

\keywords{Reaction--Diffusion, Travelling waves, tumour growth}

\maketitle

\section{Introduction}

Systems of reaction-diffusion partial differential equations (PDEs) have been used to model the interaction of species in ecology, biology and chemistry since the pioneering work of \citet{fisher1937wave} and \citet{kolmogorov1937etude}.  The canonical example is the Fisher--KPP equation, which couples linear diffusion with logistic growth, and may be written (nondimensionally) as 
\begin{equation}
\label{eq:fkpp}
\pfrac{u}{t} = \pfrac{^2u}{x^2} + u(1-u), \qquad t > 0.
\end{equation}
Generally speaking, $u$ could model a population of cells, chemical concentration, or the prevalence of a gene, for example, and the problem is considered on an infinite spatial domain ($-\infty < x < \infty$), or a semi-infinite domain ($0 < x < \infty$) with suitable (usually no-flux) boundary condition.
  
The most widely-studied property of the Fisher--KPP equation \eqref{eq:fkpp} is that a localised initial condition on an unbounded spatial domain will evolve to a travelling wave solution, by which the solution transitions from the unstable state $u=0$ to the stable state $u=1$.  It is well known that, for an exponentially decaying initial condition $u \sim \e^{-ax}$, $x\to\infty$, the travelling wave (with associated wave speed $c$) that is selected is the one that has the same exponential decay rate, if such a solution exists, or the one with maximum decay rate $a=1$ otherwise \citep{aronson1975nonlinear,aronson1978multidimensional}.  This selection implies that the wave speed is given by
\[
c = \begin{cases} 2, & a \geq 1, \\\displaystyle a + \frac{1}{a}, & a < 1. \end{cases}
\]
A compactly supported initial condition (which may be thought of as the limiting case $a\to\infty$) tends to the travelling wave with minimum speed $c=2$.  
It was established in \citet{bramson1983convergence} that, for a compactly supported initial condition, the rate at which the travelling wave speed is approached is logarithmic; if $x_f(t)$ is the spatial location of a characteristic point on the wave (for example, the point at which $u = 1/2$), then
\begin{equation}
x_f = 2t - \frac{3}{2}\log t + O(1), \qquad t \to \infty.
\end{equation}
(see also \citet{hamel2013short,nolen2017convergence}).  As well as of inherent interest, the knowledge that convergence can be of this form is useful for estimating the long-time wave speed from numerical simulations.

Many extensions to the Fisher--KPP equation \eqref{eq:fkpp} have been considered, with particular focus on the resulting wave speed of travelling wave solutions.  A modification of the reaction term in \eqref{eq:fkpp} to have a small `cut-off' was considered in \citet{brunet1997shift}:
\begin{equation}
\pfrac{u}{t} = \pfrac{^2u}{x^2} + u(1-u)f(u;\epsilon), \qquad f(u;\epsilon) = \begin{cases} 1 & u > \epsilon \\ 0 & u \leq \epsilon \end{cases},
\label{eq:fkppCutoff}
\end{equation}
where $0 < \epsilon \ll 1$.  The motivation behind this cut-off term comes from considering \eqref{eq:fkpp} as the continuum limit of a discrete process, where $\epsilon$ represents the discrete cell size. The matched asymptotic analysis in \cite{brunet1997shift} established that the correction to the minimum wave speed (that is, the speed attained by initial conditions with sufficiently fast decay) in the limit $\epsilon \to 0^+$ was logarithmic:
\begin{equation}
c \sim 2 - \frac{\pi^2}{(\log\epsilon)^2}, \qquad \epsilon \to 0^+,
\label{eq:fkppCutoffSpeed}
\end{equation}
and that this correction was universal for a more general class of cut-off functions.  The matched asymptotic argument has more recently been made rigorous using geometric singular perturbation theory in \citet{dumortier2007critical}.

The growth of malignant tumours has long been connected with the tendency of tumour cells to acidify their local extra-cellular environment \citep{boedtkjer2020acidic,stubbs2000causes}.  This tendency is due to the Warburg effect, by which cancer cells preferentially obtain energy by lactic fermentation, even when oxygen is available.  The resulting environment is unfavourable to healthy cells and leads to degradation of the extracellular matrix, resulting in more space for invading cancer cells to invade.  

Multiple approaches have been used to mathematically model the process of acid-mediated tumour growth.  Most relevant for this study is the use of reaction-diffusion partial differential equation (PDE) systems, as originated in \citet{gatenby1996reaction}.  These are generally macroscopic models in which cell populations and extracellular chemical concentrations are treated as continuous functions of space, and heterogeneities due to tissue structure or vascularity are neglected.
The modelling approach of \citet{gatenby1996reaction} has been subsequently extended in a number of ways, for example to include direct interactions between cancer and healthy cells (that is, separating the dynamics of healthy cells from degradation of extracellular matrix) \citep{martin2010tumour,mcgillen2014general}, or the effects of chemotaxis (repulsion of cells from areas of high acidity) \citep{eckardt2024mathematical}.
Multiscale models that include the coupling of extra- and intracellular acidity, often with stochastic effects included, have been developed, simulated and long-time existence proven \citep{hiremath2017mathematical,hiremath2016stochastic,hiremath2018coupled}.

The reaction-diffusion PDE model of \citet{gatenby1996reaction} may be written in nondimensional form as
\begin{subequations}
\label{eqs:gg}
\begin{align}
\pfrac{u}{t} &= \pfrac{}{x}\left[(1-v)\pfrac{u}{x}\right] + u(1-u), \\
\pfrac{v}{t} &= \gamma_1 v(1-v-\alpha h), \\
D \pfrac{h}{t} &= \gamma_2(u-h) + \pfrac{^2h}{x^2},
\end{align}
\end{subequations}
where $u(x,t)$ is the invading (cancer) cell density, $v(x,t)$ the resident (healthy) cell density, and $h(x,t)$ the acid concentration, while $\alpha$, $\gamma_1$, $\gamma_2$, and $D$ are nondimensional parameters.  In this model the resident cell population $v$, which effectively represents the combination of healthy cells and extracellular matrix, is stationary (non-diffusive) and proliferates with a carrying capacity that decreases with the acid concentration.  The acid $h$ is produced by the invading cells, while undergoing diffusion and degradation.  The parameter $D$ represents the ratio of cell to acid diffusivities, and would generally be taken to be small.  Importantly, the diffusivity of the invading cell population is a decreasing function of the resident cell population, such that when the healthy cell population is at capacity, the diffusivity vanishes.  The decrease in resident cell population due to acidity is thus required to allow the tumour to grow.  

In \citet{gatenby1996reaction} it was numerically observed that solutions to their three-component model \eqref{eqs:gg} tend to travelling wave solutions as time increases.  In such solutions, there is an advancing front of the invading cells $u$ similar to the Fisher--KPP model coupled to a receding front of the resident cells $v$.  A particular point of interest was that when the death of resident cells due to acidity is large (large $\alpha$ in \eqref{eqs:gg}), there is a separation between the advancing and receding fronts, in between which both cell populations are small.  Such a region corresponds to a phenomenon sometimes observed in biological experiments referred to as the hypocellular \textit{interstitial gap}.  \citet{gatenby1996reaction} estimated that the size of this region was logarithmically large in $\alpha$.  The asymptotic analysis of travelling wave solutions has also been carried out more recently in \citet{fasano2009slow} and \citet{davis2022traveling}; the interstitial gap is a feature of waves that have speed $O(1)$ in the limit that acid diffusivity is much larger than the cell diffusivity.  These waves are referred to as `slow' waves in \citet{fasano2009slow} \citet{davis2022traveling}, the `fast' waves relating to the finite diffusivity ($D > 0$ in \eqref{eqs:gg}) of the acid species.

Simplified versions of the three-species model of \citet{gatenby1996reaction} have been studied, by neglecting the acid species, and directly modelling the effect of resident cell population on the diffusivity of the invading species.  \citet{gallay2022propagation} and \citet{mascia2024numerical} consider existence of travelling wave solutions for the model
\begin{equation}
\label{eq:mascia}
\pfrac{u}{t} = \pfrac{}{x}\left[(1-v)\pfrac{u}{x}\right] + u(1-u), \qquad \pfrac{v}{t} = \gamma_1 v(1 - v - \alpha u).
\end{equation}
Formally this model is recovered from \eqref{eqs:gg} by taking the coupling between invading cell population $u$ and acid concentration $h$ to be strong, that is, $\gamma_2 \gg 1$, so that $u \sim h$.  The stability of travelling wave solutions of this system is considered in \citet{swartwood2025stability}, who show that all travelling waves are stable, regardless of the speed $c$; this is similar to the Fisher--KPP equation \eqref{eq:fkpp}, where the marginal stability threshold is considered as the relevant wavespeed selection mechanism \citep{van1988front,van2003front}.  

\citet{colson2021travelling} consider a version of the reduced model \eqref{eq:mascia} in which the resident cell species does not recover:
\begin{equation}
\label{eq:colson}
\pfrac{u}{t} = \pfrac{}{x}\left[(1-v)\pfrac{u}{x}\right] + u(1-u), \qquad \pfrac{v}{t} = -\gamma uv.
\end{equation}
This may be considered the large-$\alpha$ approximation of \eqref{eq:mascia}, identifying $\gamma = \alpha \gamma_1$.  In this model the initial resident cell population plays an important role, as the initial population may or may not be taken to be the value at which the diffusivity of the invading cell population vanishes ($v=1$, nondimensionally).  In \citet{colson2021travelling} families of travelling wave solutions are constructed in phase space by a shooting algorithm.  A minimum wave speed is identified as a function of the resident cell death rate, and the initial resident cell population. 

Recently, \citet{browning2019bayesian} and \citet{el2021travelling} consider a similar model to \eqref{eq:colson}, in which both the diffusion and the carrying capacity of the invading cells is inhibited by the healthy cell population:
\begin{subequations}
\label{eqs:PDEs}
\begin{align}
\frac{\partial u}{\partial t} &= \frac{\partial}{\partial x}\left[(1-v)\frac{\partial u}{\partial x} \right] + u(1-u-v), \label{eq:uPDE} \\
\frac{\partial v}{\partial t} &= -\gamma u v \label{eq:vPDE},
\end{align}
\end{subequations}
with $t>0$.  In the context of travelling waves this model is posed on an infinite domain $-\infty < x < \infty$, or equivalently a semi-infinite domain $0 < x < \infty$ with appropriate no-flux condition at $x=0$.  This model may be seen as a version of the system \eqref{eq:colson} with an additional effect of resident cell density $v$ on the carrying capacity of invading cells $u$.  Nominally this model has a single parameter, the dimensionless resident cell death rate $\gamma$, but since in this model there is no proliferation of $v$, the initial condition of $v$ again plays an important role.  The diffusivity of $u$ vanishes at $v=1$, but the diffusion will only be degenerate if the initial condition takes this value.  Although travelling wave solutions of \eqref{eqs:PDEs} were observed in PDE simulations in \citet{el2021travelling} over a range of parameter values and initial conditions, and the boundary value problem for travelling wave solutions was posed, the explicit computation of travelling wave solutions was not performed, and the mechanism by which the wave speed is determined by the parameters was not established.

In this article we focus on the model \eqref{eqs:PDEs}, with the aim to build insight into not only the properties of travelling wave solutions, but also the selection of a given wave speed from a given initial condition.  While this was only explored using PDE simulations in \citet{el2021travelling}, we explicitly calculate travelling wave solutions via numerical continuation over a large range of parameter values (in particular the resident cell death rate $\gamma$), and compare with time-dependent numerical simulations.  We are thus able to observe that, although the details of the phase space are more complicated, the principle by which the wave speed is selected is the same as for the simpler Fisher--KPP system: that is, an initial condition with sufficiently small exponential decay rate will tend to the travelling wave solution with the same decay rate (if it exists), and that an initial condition with compact support (or sufficiently large exponential decay) will tend to the travelling wave solution that has maximum decay rate.  This is similar to the observation made in \citet{colson2021travelling} for their simpler system \eqref{eq:colson}.

In addition, we explore the asymptotic limit as the resident cell death rate $\gamma$ becomes large, which was not achieved in \citet{el2021travelling}.  As discussed above, it is in this limit that we may connect the hierarchy of models \eqref{eqs:PDEs} and \eqref{eq:colson} to \eqref{eq:mascia} and ultimately the complete Gatenby--Gawlinski model \eqref{eqs:gg}.  In addition, this limit corresponds to the large-$\alpha$ limit of the Gatenby--Gawlinski model in which we expect to see a separation of advancing and receding fronts.  Importantly, in this limit each of the above models becomes a perturbed version of the Fisher--KPP equation \eqref{eq:fkpp}, in that the equation for invading cells $u$ reduces to \eqref{eq:fkpp} when the resident population $v$ vanishes.

We demonstrate that there is a very strong connection between the asymptotic limit of large resident cell death rate and the analysis of the Fisher--KPP system with cut off \eqref{eq:fkppCutoff} \citep{brunet1997shift, dumortier2007critical}, with $\gamma^{-1}$ playing the role of the small parameter $\epsilon$.  The leading order correction to the velocity is the same as the universal correction \eqref{eq:fkppCutoffSpeed}.  We also find the next term in this relation, although this requires the numerical calculation of prefactors in each of the leading order problems in the asymptotic analysis.  This asymptotic solution also exhibits a region between fronts where both populations are small, that is, an interstitial gap; the width of this zone is also logarithmically dependent on the cell death parameter.  As our numerical continuation method allows us to numerically compute travelling wave solutions up to very large ($10^{11}$) values of $\gamma$, we are able to numerically observe the correctness of the asymptotic approximation.  The universality of this asymptotic approach demonstrates the fundamental dynamics that is at play in interstitial gap formation in the large family of reaction-diffusion models of acid-mediated cancer invasion.

\section{Model and numerical PDE simulations}
\label{sec:modelAndNumerics}

We start by describing the method and results of numerical simulations of the system \eqref{eqs:PDEs}.  \citet{el2021travelling} report on a large number of numerical simulations of this system; our aim in this section is to reproduce their main observations in order to confirm the validity of the travelling wave solutions in the next section.  The system \eqref{eqs:PDEs} comprises a population density of invading cells $u(x,t)$ invading a stationary population of resident cells $v(x,t)$, where both the diffusive flux and carrying capacity of the invading cells is negatively affected by the resident population.  The nondimensional parameter $\gamma$ in \eqref{eqs:PDEs} represents the rate of destruction of resident cells, compared to the  time scale of invading cell proliferation.  \citet{browning2019bayesian} and \citet{el2021travelling} contain further details of this model, including the nondimensionalisation.

Initially, the domain consists of a spatially localised invading cell population $u_0(x)$ in a uniform resident cell population $v_0$:
\begin{equation}
u(x,0) = u_0(x), \qquad v(x,0) = v_0. \label{eq:IC}
\end{equation}
Note that if $v_0 = 0$ identically, the system \eqref{eqs:PDEs} reduces to the Fisher--KPP problem for $u$.  In addition, the diffusivity vanishes when $v=1$; in this paper, we will consider only $0 < v_0 < 1$, as our main focus is on the effects of $v$ on the proliferation, rather than the effect of degenerate diffusion (see the discussion in Sect.~\ref{sec:discussion}).

To simulate \eqref{eqs:PDEs} numerically, we discretise in space using a standard cell-centred finite volume method \citep{patankar1980numerical} with equally sized cells, and advance in time using the \texttt{ode15s} algorithm in MATLAB.  To observe the formation of travelling wave behaviour, we consider a large spatial domain $x \in [0, L]$ with $L \gg 1$ approximating a semi-infinite domain, with zero flux conditions ($\partial_x u = 0$) conditions imposed at the boundaries.  A nonzero initial condition $u(x,0)$ near $x=0$ is then used to represent a localised initial condition.  We generally choose $L$ to be between $150$ and $200$ and $N = 10^4$ grid points, which is more than adequate for the numerical simulations to have converged.

We calculate the wave speed from PDE simulations by estimating the location of the travelling wave front $x_f(t)$ where $u(x_f(t),t) = 0.5$, interpolating the solution between neighbouring grid points if required.  This calculation is sensitive in that sufficient time must have passed to minimise the effect of the initial condition, and the domain size $L$ must be sufficiently large to minimise the effect of the right hand boundary.  Rather than fit a simple linear function to $x_f$ to estimate the wave speed \citep[e.g.]{el2021travelling}, we follow the work of \citet{bramson1983convergence} on the Fisher--KPP equation (see also \citet{nolen2017convergence}), who showed that the wave speed is approached only logarithmically as $t\to \infty$.  Correspondingly, we expect that a constant wave speed is approached according to
\begin{equation}
x_f(t) \sim ct + k_0\log(t) + k_1, \qquad t \to \infty,
\label{eq:waveSpeedFit}
\end{equation}
where $c$ is the wave speed, and $k_0$ and $k_1$ are constants; $k_1$ depends on the initial condition, and while $k_0 = -3/2$ for the minimum-speed Fisher--KPP solution, it likely depends on the parameters for the more complex model \eqref{eqs:PDEs}.  To estimate $c$ we thus fit the nonlinear relationship \eqref{eq:waveSpeedFit} for large times using a least-squares procedure.

\begin{figure}
\centering
\includegraphics[width=\textwidth]{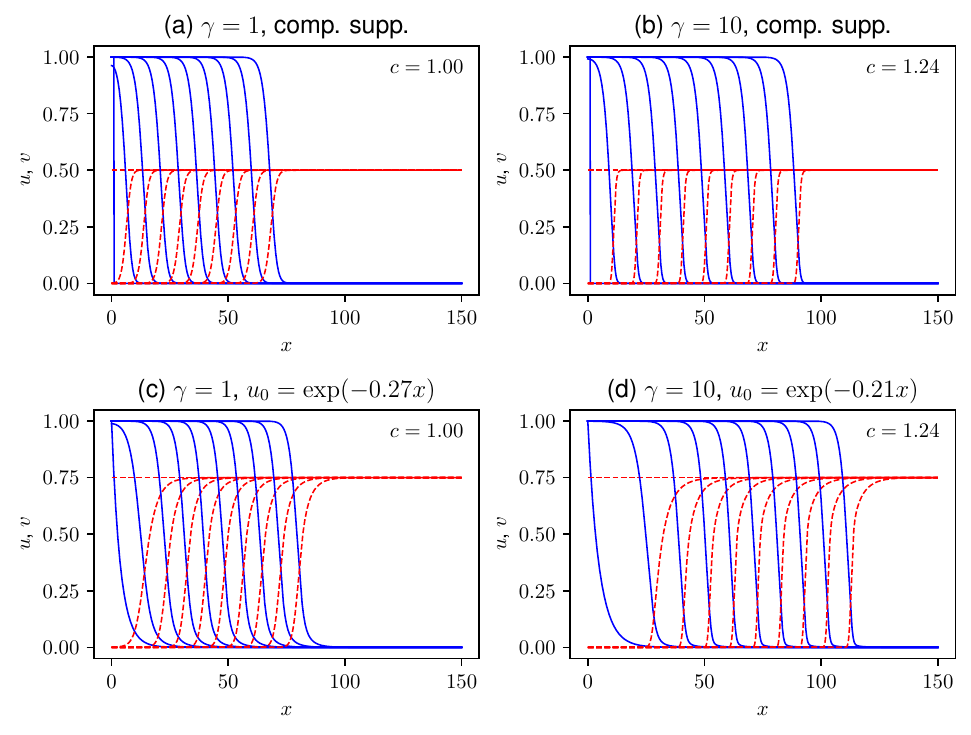}
\caption{Numerical solutions of the PDE system \eqref{eqs:PDEs}, with various initial conditions and parameter values.  In (a, b), $v_0 = 0.5$ and $u_0$ is compactly supported, while $\gamma = 1,10$, respectively.  In (c, d), $v=0.75$ and $u_0(x) = \e^{-ax}$ is an exponentially decreasing initial condition, where the value of $a$ has been chosen specifically to result in the same wave speed as for the results in (a) and (b), respectively.  The estimated wave speed $c$ (from \eqref{eq:waveSpeedFit}) is noted in each simulation.}
\label{fig:PDESolutions}
\end{figure}

In Fig.~\ref{fig:PDESolutions} we plot typical results (profiles of $u$ and $v$ over time) for a combination of different cell death rates $\gamma$ and initial conditions.  We choose specific combinations of $\gamma$, $v_0$, and either compactly supported or exponentially decaying $u_0(x)$, in order to demonstrate the nontrivial dependence of the wave speed on these parameters that will be further explored using the travelling wave analysis in the next section.

In Fig.~\ref{fig:PDESolutions}(a,b), we consider an initially compactly supported invading cell population:
\[
u_0(x) = \begin{cases} 1 & x < 1 \\ 0 & x\geq 1\end{cases}, \qquad v_0 = 0.5,
\]
for the values $\gamma = 1$ and $\gamma=10$, respectively.  It is apparent from examining the behaviour of solutions that the system is tending toward a travelling wave of fixed shape and constant speed.  The speed $c$ depends on the value of $\gamma$; applying the fitting procedure described above to the function \eqref{eq:waveSpeedFit}, we estimate $c=1.00$ for $\gamma=1$, and $c=1.24$ for $\gamma = 10$.  We note that these solutions have wave speed considerably less than the minimum wave speed $c=2$ of the Fisher--KPP system \eqref{eq:fkpp}, although $c$ is seen to increase as $\gamma$ increases.

In Figs.~\ref{fig:PDESolutions}(c,d) we consider initial conditions with exponential decay, and larger initial resident cell population $v_0$:
\[
u_0(x) = \e^{-ax}, \qquad v_0 = 0.75,
\]
for $a = 0.27$ and $a = 0.21$, respectively, and again for $\gamma = 1$ and $\gamma = 10$, respectively.  These solutions are seen to converge to travelling waves with the same respective wave speeds: $c=1$, and $c=1.24$, respectively.  These results show that different combinations of initial conditions can lead to the same wave speed.  We will further explain this relationship in the travelling wave analysis of the next section, using these parameter values as specific examples.

%
%
%
%
%

\begin{figure}
\centering
\includegraphics[width=\textwidth]{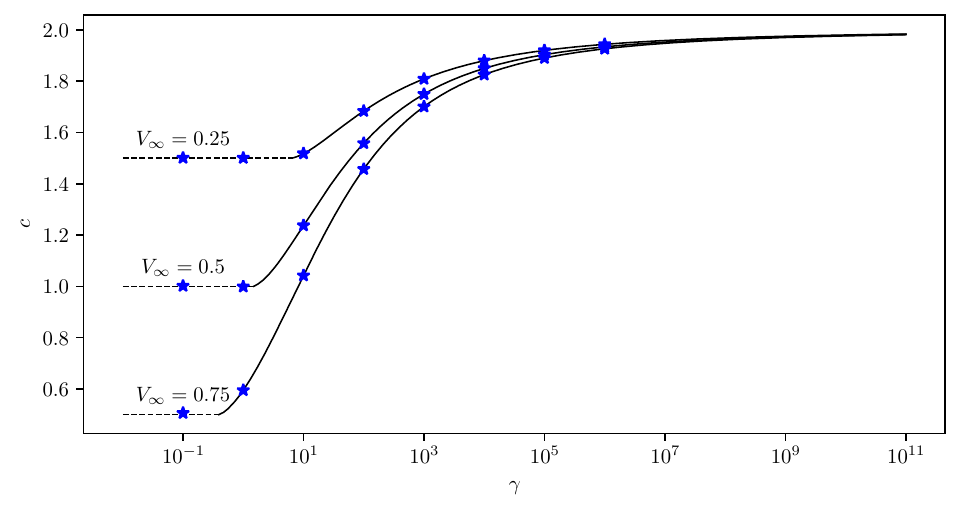}
\caption{The dependence of wave speed $c$ on death rate $\gamma$, for different values of initial resident cell population $v_0 = V_\infty$.  Blue symbols indicate the estimation of wave speed from PDE simulations with compactly support initial conditions, described in Sect.~\ref{sec:modelAndNumerics}.  The dashed curves represent the predicted wave speed according to the condition that $V_\infty = V_c$ \eqref{eq:Vc}, while the solid curves represent the predicted wave speed according to the condition $V_\infty = V_s$, where $V_s$ must be found by numerically solving the system \eqref{eqs:odeSystem}.  The PDE simulations with initial condition $v_0 = V_\infty$ predict the same wave speed as the travelling wave analysis.}
\label{fig:waveSpeedComparison}
\end{figure}

To further demonstrate the possible behaviour of the system \eqref{eqs:PDEs}, we run simulations and estimate travelling wave speeds for a much wider range of parameter values and initial conditions.  For compactly supported $u_0(x)$, we estimate the wave speed for values of $\gamma$ ranging from $0.1$ to $10^6$, and for values $v_0 \in \{ 0.25, 0.5, 0.75 \}$.  These estimated wave speeds are plotted in Fig.~\ref{fig:waveSpeedComparison}, and will be used to validate the predictions made from the direct computation of travelling wave solutions in the next section.  The most notable property of these speeds is that for sufficiently small $\gamma$, the wave speed is  independent of $\gamma$, while for sufficiently large $\gamma$, the wave speed becomes dependent on $\gamma$.  The point of transition between these behaviours depends on $v_0$, but in all cases, as $\gamma\to\infty$, the speed $c$ approaches $2$ from below.

\begin{figure}
\centering
\includegraphics[width=\textwidth]{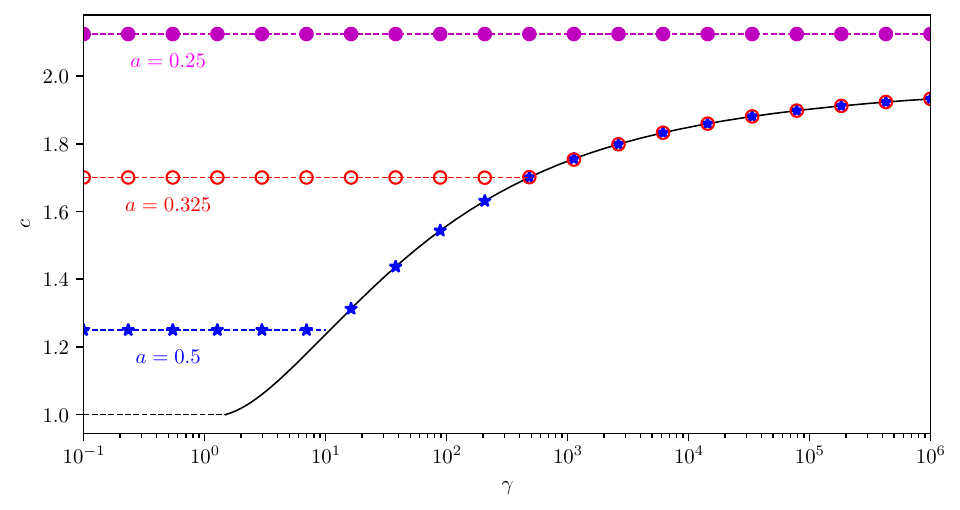}
\caption{The dependence of wave speed $c$ on death rate $\gamma$, for initial resident cell population $v_0 = 0.5$ and different invading cell initial conditions $u_0 = \exp(-ax)$, for $a \in \{0.25, 0.325, 0.5\}$.  The symbols represent numerical solutions for $a=0.25$ (filled circles), $a=0.325$ (unfilled circles), and $a=0.5$ (stars), respectively.  The wave speed predicted by the formula \eqref{eq:speed2} for each $a$ is depicted by the coloured dashed lines.  The wave speed of a compactly supported initial condition, corresponding to that in Fig.~\ref{fig:waveSpeedComparison}, is shown as black dashed/solid lines.}
\label{fig:waveSpeedvsIC}
\end{figure}

Finally, we further test the effect of initial condition of $u$ by calculating the wave speed for $\gamma \in [0.1, 10^6]$, fixed $v_0 = 0.5$, and exponentially decaying initial conditions $u_0 = \e^{-ax}$ for $a \in \{0.5, 0.325, 0.25\}$.  These wave speeds are plotted in Fig.~\ref{fig:waveSpeedvsIC}.  These results show that for sufficiently small $a$ and $\gamma$, the wave speed is again independent of $\gamma$.  Furthermore, for a fixed $a$ sufficiently small such that the wave speed $c$ is greater than 2, for example $a=0.25$ shown in Fig.~\ref{fig:waveSpeedvsIC}, the wave speed is independent of $\gamma$ over all $\gamma$.  On the other hand, if $a$ is sufficiently large such that the speed is less than $2$, for sufficiently large $\gamma$ the wave speed is the same as the $\gamma$-dependent wave speed of the compactly supported initial condition.  Equivalently, for a fixed $\gamma$, as $a$ is increased, the wave speed is ultimately limited by the corresponding wave speed of the compactly supported initial condition of the same value of $\gamma$.   These phenomena will also be explained by the travelling wave analysis in the next section.

\section{Calculation of travelling wave solutions}
\label{sec:travWave}

In this section we will compute travelling wave solutions to \eqref{eqs:PDEs}, and analyse their far field behaviour that determines which travelling wave solutions are selected from a given initial condition.    Let $z = x-ct$ be the travelling wave coordinate, and let $u(x,t) = U(z)$, $v(x,t) = V(z)$, and $\partial _x u(x,t) = W(z)$.  The system of PDEs then becomes the system of three autonomous differential equations:
\begin{subequations}
\label{eqs:odeSystem}
\begin{align}
\frac{\mathrm dU}{\mathrm dz} &= W, \label{eq:Ueq}\\
\frac{\mathrm dV}{\mathrm dz} &= \frac{\gamma }{c} UV, \label{eq:Veq} \\
\frac{\mathrm dW}{\mathrm dz} &= \frac{1}{1-V}\left(\frac{\gamma}{c} UVW - cW - U(1-U-V) \right). \label{eq:Weq}
\end{align}
\end{subequations}
This system \eqref{eqs:odeSystem} has an isolated equilibrium at $(U,V,W) = (1,0,0)$.  In addition, each point on the $V$-axis $(U,V,W) = (0,V_\infty,0)$ is an equilbrium, each of which is therefore not isolated.  Travelling wave solutions correspond to heteroclinic orbits connecting $(1,0,0)$ (as $z\to-\infty$) to a point on the $V$-axis (as $z\to\infty)$.  The travelling wave solution relevant to an initial value problem to \eqref{eqs:PDEs} with $v(x,0) = v_0$ is one in which this point on the $V$-axis has $V_\infty=v_0$.

The eigenvalues and eigenvectors of the fixed point $(1,0,0)$ are
\begin{subequations}
\begin{equation}
\lambda_1 = \frac{\gamma}{c}, \quad \lambda_2 = \frac{-c + \sqrt{c^2 + 4}}{2}, \quad \lambda_3 = \frac{-c - \sqrt{c^2 + 4}}{2}.
\end{equation}
\begin{equation}
\bm E_1 = \begin{bmatrix} 0 \\ 1 \\ 0 \end{bmatrix}, \quad \bm E_{2,3} = \begin{bmatrix} 1 \\ 0 \\ \lambda_{2,3} \end{bmatrix}.
\end{equation}
\label{eqs:eigsAtNegativeInfinity}
\end{subequations}
This equilibrium thus has a two-dimensional unstable manifold.  At the other end, the eigenvalues near a point $(0,V_\infty,0)$ are
\begin{equation}
\lambda_1' = 0, \quad \lambda_2' = -\mathcal C + \sqrt{\mathcal C^2 - 1}, \quad \lambda_3' = -\mathcal C - \sqrt{\mathcal C^2 - 1}, \quad \mathcal C = \frac{c}{2(1-V_\infty)},
\label{eq:eigsAtInfinity}
\end{equation}
with eigenvectors
\[
\bm E_1' = \begin{bmatrix} 0 \\ 1 \\ 0 \end{bmatrix}, \quad \bm E'_{2,3} = \begin{bmatrix} \lambda'_{2,3} \\ \gamma V_\infty/ c \\ (\lambda'_{2,3})^2 \end{bmatrix}.
\]
The zero eigenvalue $\lambda_1'$ (with eigenvector pointing along the $V$-axis) is due to the non-isolated nature of each fixed point.  The other two eigenvalues are negative, indicating that trajectories are attracted to the $V$-axis in its neighbourhood.
Because the unstable manifold is two-dimensional, there is one degree of freedom in trajectories leaving $(1,0,0)$.  There is thus a one-parameter family of trajectories that connect $(1,0,0)$ to the $V$-axis, each ending at a different value of $V_\infty$ between $0$ and $1$.  


We compute solutions of \eqref{eqs:odeSystem} numerically by solving it as a boundary value problem on a large domain $z \in (0, z_\infty)$, using the numerical continuation package \textsc{Auto}-07p \citep{doedel2007auto}.  The use of numerical continuation will allow us to compute solutions for extreme values of parameters (in particular, very large $\gamma$).  Appropriate boundary conditions are
\[
U(0) = 1 + \varepsilon, \quad W(0) = -\varepsilon\lambda_2, \quad V(z_\infty) = V_\infty,
\]
where $\varepsilon$ is a small parameter and $z_\infty \gg 1$; in our calculations we take $\varepsilon = 10^{-5}$ and $z_\infty = 50$.  These boundary conditions enforce that near $(1,0,0)$, trajectories are on the eigenvector with the smaller positive eigenvalue.  For a given $\gamma$ and $c$ we are able to construct travelling wave solutions for any value of $V_\infty$ between 0 and 1, by starting at $V_\infty=0$ (at which point $V$ vanishes identically, and $U,W$ are solutions to the Fisher--KPP system), and then increasing $V_\infty$ by numerical continuation.  

\begin{figure}
\centering
\includegraphics[width=\textwidth]{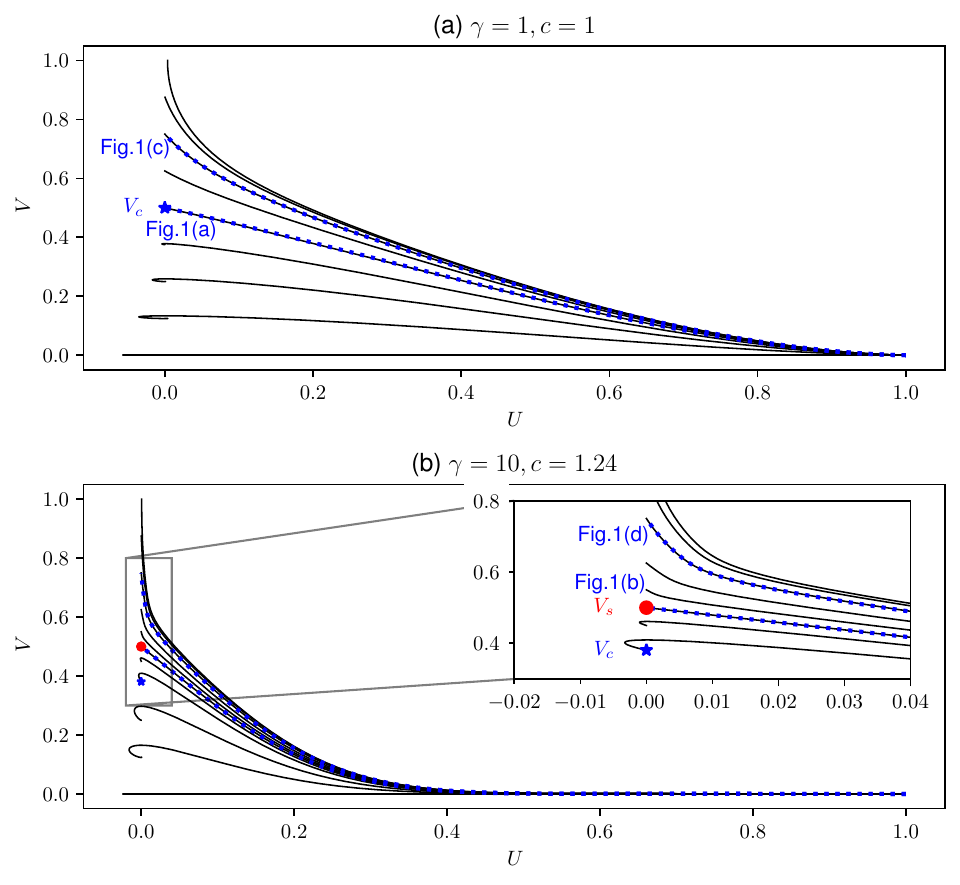}

\caption{Trajectories in $(U,V,W)$ phase space, projected onto the $(U,V)$ plane, for (a) $(\gamma, c) = (1,1)$ and (b) $(\gamma, c) = (10, 1.24)$.  In (a), the parameters are such that trajectories spiral into the $V$-axis for $V < V_c$, and monotonically approach the $V$-axis for $V \geq V_c$.  In (b), there is also a critical value $V_s$ such that $U$ approaches zero from below when $V < V_s$ and above when $V > V_s$.
Dotted curves marked (a)--(d), corresponding to the late-time PDE simulations depicted in Figs.~\ref{fig:PDESolutions}a--d, respectively, show that PDE simulations are attracted to travelling wave solutions with either the same decay rate as initial condition (if such a solution exists), or to the travelling wave with fastest decay rate, at $V_\infty = \max\{V_c, V_s\}$.
}
\label{fig:phasePlanes}
\end{figure}

We plot examples of such trajectories, projected onto the $(U,V)$ plane, in Fig.~\ref{fig:phasePlanes}.  In this figure we choose two sets of parameters $(\gamma, c) = (1,1)$ and $(10, 1.24)$ in order to compare to the travelling wave solutions observed in the PDE simulations shown in Fig.~\ref{fig:PDESolutions}.  
For given values of $\gamma$ and $c$, there is a one-parameter family of travelling wave solutions, that we may consider as characterised by the value $V_\infty$.  Since the wave speed $c$ is a parameter of the system, we thus expect that for a given $V_\infty$ there is a one-parameter family of travelling wave solutions, parametrised by $c$.  In Fig.~\ref{fig:phasePlanes} we observe that the travelling waves that emerge from the PDE simulations do indeed correspond to certain trajectories in phase space.  However, we are yet to determine which trajectory (in particular, which wave speed $c$) will correspond to the long-time behaviour of a solution to the PDE system \eqref{eqs:PDEs} for a given initial condition.
In order to determine this wave-speed selection, we need to consider the behaviour of trajectories near the front of the travelling wave, that is, near the $V$-axis.   

For initial conditions of the form $u(x,0) \sim \e^{-ax}$, $x\to\infty$, as shown in the previous section, one may expect the appropriate travelling wave to be that which has the same exponential decay.  Since this decay is given by the least negative eigenvalue $\lambda_2'$ in \eqref{eq:eigsAtInfinity}, by equating $a = -\lambda_2'$ we find the relationship
\begin{equation}
c = (1-V_\infty)\left(a + \frac{1}{a}\right).
\label{eq:speed1}
\end{equation}
For initial conditions that have sufficiently slow decay this relationship does indeed hold.  For example, the solutions shown in Fig.~\ref{fig:PDESolutions}(c) and Fig.~\ref{fig:PDESolutions}(d), corresponding to values $(v_0, \gamma, c) = (0.75, 1, 1)$, and $(v_0, \gamma, c) = (0.75, 10, 1.24)$, satisfy this relationship.  In the phase plane in Fig.~\ref{fig:phasePlanes} it is also observed that these PDE simulations correspond to trajectories in phase space that terminate at $V_\infty = 0.75$ for the appropriate values of $\gamma$ and $c$.  This relationship can also be seen in the results of simulations depicted in Fig.~\ref{fig:waveSpeedvsIC}; for values of $a$ such that $c > 2$, that is, 
\begin{equation}
a < \frac{1}{1-V_\infty} - \sqrt{\frac{1}{(1-V_\infty)^2} - 1},
\label{eq:stupidSpeed}
\end{equation}
the wave speed is always determined by the criterion \eqref{eq:speed1}.  For larger $a$, this relation holds for sufficiently small $\gamma$, until a greater constraint on the wave speed is encountered.

As occurs in the analysis of the Fisher--KPP system, a further criterion results from avoiding parameters that result in travelling waves in which $U$ does not remain positive over all $z$.  This will be important for sufficiently large decay rate $a$, and in particular compactly supported initial conditions.  Firstly, we note that for a given wavespeed $c$, there is a critical value $V = V_c$, at which the eigenvalues $\lambda'_{2,3}$ switch from complex to real:
\begin{equation}
V_c = 1 - \frac{c}{2}.
\label{eq:Vc}
\end{equation}
This value corresponds to $\mathcal C = 1$ in \eqref{eq:eigsAtInfinity}. Notably, $V_c$ is independent of the decay rate $\gamma$.  Trajectories for which $V_\infty < V_c$ spiral into the relevant point on the $V$-axis, and since $U$ will become negative on such a spiral, these trajectories cannot correspond to solutions to \eqref{eqs:PDEs} that have evolved from a non-negative initial condition.  This property is determined locally, and is equivalent to the condition used to determine the minimum wave speed in the Fisher--KPP equation.  It results in the prediction for the wave speed:
\begin{equation}
c = 2(1-V_\infty), \qquad V_\infty = V_c.
\label{eq:speed2}
\end{equation}
This relation can also be observed in the PDE simulations; it is the relevant constraint for the solution $v_0 = 0.5, \gamma = 1, c=1$ as depicted in Fig.~\ref{fig:PDESolutions}(a), and as depicted by the dashed lines in Fig.~\ref{fig:waveSpeedComparison}, for sufficiently small $\gamma$.  

However, \eqref{eq:speed2} does not explain the dependence of wave speed on $\gamma$ when $\gamma$ is sufficiently large, so the wave speed in this regime must be limited by a different criterion.  The reason for this additional constraint is that for some parameter regimes, $U$ becomes negative as $z\to\infty$, even when $V_\infty > V_c$, that is, when the eigenvalues of the steady state are real.  In Fig.~\ref{fig:phasePlanes}(b) (in particular the inset of that figure), in which $c = 1.24$ and $\gamma=10$, we observe that trajectories for which $V_\infty$ is slightly larger than $V_c$ still approach the $V$-axis from the negative $(U<0)$ side, despite the fact that eigenvalues are real at these points.  In this parameter regime there is a point $V_s$ (indicated as a red dot in Fig.~\ref{fig:phasePlanes}(b)) that separates trajectories that have $U$ become negative, from those that remain positive.  Also seen in this figure is that this is the appropriate trajectory corresponding to the PDE simulation depicted in Fig.~\ref{fig:PDESolutions}(b), where $(v_0, \gamma, c) = (0.5, 10, 1.24)$, which does not satisfy the condition \eqref{eq:speed2}.  The existence of such a threshold value was also noted in the similar model of \citet{colson2021travelling} \eqref{eq:colson} described in the introduction.  Unlike the point $V_c$, the point $V_s$ is not determined as a local property of the phase space (it is a property of the manner in which the unstable manifold of $(1,0,0)$ intersects the $V$-axis), so is unlikely to have a closed-form expression.

We may understand the existence of this critical point $V_s$ by considering the linearisation near the $V$-axis. For $V_c < V_\infty < 1$, a trajectory near the $V$-axis will behave as
\begin{equation}
\begin{bmatrix} U \\ V \\W \end{bmatrix} \sim \begin{bmatrix}0 \\ V_\infty \\ 0 \end{bmatrix} + r_2 \bm E_2'\e^{\lambda_2' z} + r_3 \bm E_3'\e^{\lambda_3' z}, \qquad z\to\infty,
\label{eq:nearVaxisExpansion}
\end{equation}
where $\lambda_{2,3}'$ and $\bm E_{2,3}'$ are the eigenvalues and eigenvectors \eqref{eq:eigsAtInfinity}.  For $V>V_c$ , $\lambda_3' < \lambda_2' < 0$, so that the generic case is that $U\to 0$ with decay rate given by the less negative eigenvector $\lambda_2'$.
The magnitudes of one of the constants $r_2$ and $r_3$ is arbitrary due to translational invariance, but these constants are otherwise unique properties of a given trajectory.  Importantly, a trajectory may have $r_2 > 0$, which means that (for the eigenvector written as in \eqref{eq:eigsAtInfinity}) the trajectory approaches the $V$-axis from the negative $U$ direction.  These are the trajectories observed in Fig.~\ref{fig:phasePlanes}(b) for $V_c < V_\infty < V_s$.  These trajectories also cannot correspond to a travelling wave solution of \eqref{eqs:PDEs} reached from a non-negative initial condition.  
The critical point $V_s$ represents the end point of the particular trajectory at which the coefficient $r_2$ in \eqref{eq:nearVaxisExpansion} goes from being positive to negative, that is, at $V_\infty = V_s$, $r_2 = 0$.  What is special about this point then is that $U\to 0$ with decay rate given by the \textit{more negative} eigenvalue $\lambda_3'$.  In other words, for parameter values in which a point $V_s$ exists, the trajectory ending at $V_s$ (rather than $V_c$) is the travelling wave solution with the greatest decay as $z\to \infty$.


We calculate the location of $V_s$ in our numerical continuation procedure by tracking the value of $U(z_\infty)$.  While in the limit $z_\infty\to\infty$ this value will vanish, for a finite $z_\infty < \infty$ it is generally small but nonzero.  In this case, the value of $V_\infty$ where $U(z_\infty) = 0$ corresponds (closely enough for the computations) to the threshold value $V_\infty = V_s$.  Using numerical continuation we are then able to fix $U(z_\infty)=0$ and allow $\gamma$ and $c$ to vary, for a given $V_\infty$.  The parameter $\gamma$ can then be increased to very large values, for which the wave speed slowly approaches $2$.  As $\gamma$ is decreased, eventually $V_s$ collides with the other critical point $V_c$.

Using this procedure we produce curves of wave speed $c$ as a function of $\gamma$, for given values of the far-field resident cell population $V_\infty$.  These results are plotted in Fig.~\ref{fig:waveSpeedComparison} along with the estimates from the PDE simulations for compactly supported initial conditions.  When $\gamma$ is sufficiently small such that $V_s$ does not exist, we instead use the condition \eqref{eq:speed2} assuming $V_\infty = V_c$.  The comparison confirms that we have correctly identified the wave speed selection mechanism.

\begin{figure}
\centering
\includegraphics[width=\textwidth]{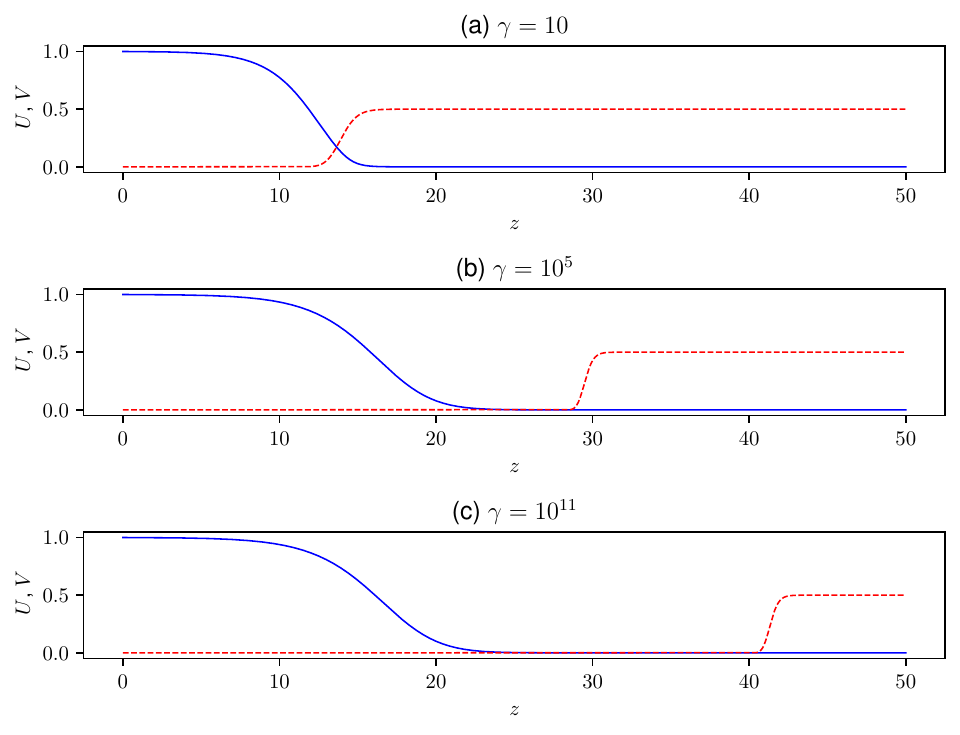}
\caption{Travelling wave solutions for $U$ (blue solid curves) and $V$ (red dashed curves), for $V_\infty = 0.5$ and (a) $\gamma = 10$, (b) $\gamma = 10^5$, and (c) $\gamma = 10^{11}$.  The interstitial gap (the region in which both cell populations are small) is seen to grow slowly; the asymptotic analysis of Section \ref{sec:asymptotics} predicts this region is of order $\log\gamma$ as $\gamma\to\infty$.}
\label{fig:bigGammaProfiles}
\end{figure}

As we can explicitly calculate travelling wave solutions, we are able to continue to much larger values of $\gamma$ than would be feasible if we relied only on PDE simulations.  We observe in Fig.~\ref{fig:waveSpeedComparison} that in the limit that $\gamma$ becomes large, the wave speed tends to a constant value $2$, but does so very slowly.  Example travelling wave profiles are shown for $V_\infty = 0.5$, and $\gamma \in\{10, 10^5, 10^{11}\}$, in Fig.~\ref{fig:bigGammaProfiles}.  These profiles highlight that as $\gamma$ increases, the travelling wave solutions exhibit a region in which both cell populations are very small, and this region grows, but again slowly, as $\gamma\to\infty$.  This region corresponds to the hypocellular interstitial gap described in \citet{gatenby1996reaction}.  In the next section we carry out the matched asymptotic analysis in this limit, which will explain the logarithmic dependence on $\gamma$ of both the wave speed and the interstitial gap width.

\section{Large death rate asymptotic analysis}
\label{sec:asymptotics}

In this section we describe the asymptotic analysis of the relation between resident cell death rate $\gamma$ and wave speed $c$ in the limit that $\gamma$ becomes large.  We restrict ourselves to the travelling wave branches that correspond to compactly supported initial conditions (that is, $V_\infty = V_s$ in the preceding section).  We introduce the wave speed correction $\delta$, such that the wave speed is $c = 2 - \delta$.  From the above numerical results, $\delta \to 0 $ as $\gamma \to \infty$, and the main task of the asymptotic analysis is to determine how $\delta$ depends on $\gamma$ in this limit.

\begin{figure}
\centering
\includegraphics[width=.85\textwidth]{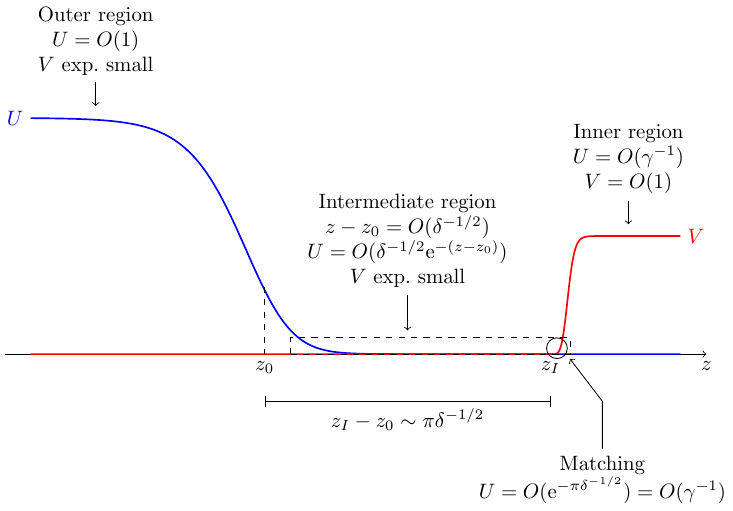}
\caption{The asymptotic structure of travelling wave solutions in the large death rate ($\gamma \to \infty$) limit.  The wave speed $c = 2-\delta$ where $\delta\to 0$ as $\gamma \to \infty$.  In the outer region, $V$ is exponentially small, and $(U,W)$ is given by the Fisher-KPP travelling wave with $c=2$.  In the intermediate region, the effect of sub-critical wave speed $(0 < \delta \ll 1)$ is felt, while $V$ remains negligible.  In the inner region, the effect of the far field condition enforces $V = O(1)$. The dependence of $\delta$ on $\gamma$ is found by matching between the intermediate and inner regions.}
\label{fig:asymptoticSchematic}
\end{figure}

The asymptotic structure that arises in the large $\gamma$ limit is depicted in Fig.~\ref{fig:asymptoticSchematic}.  The important observation is that for large $\gamma$, \eqref{eq:Veq} implies that $V$ is exponentially small except near the front (where $V\to V_\infty$ as $z\to\infty$).  There are thus three main regions: an outer region, where to leading order $U$ is given by the critical wavespeed ($c=2$) solution of the Fisher--KPP equation; an intermediate region, where $V$ is still negligible but the effect of slightly subcritical wavespeed is felt (that is, $\delta > 0$); and finally an inner region, where $V$ must be order unity to match with the far field condtion ($V\to V_\infty$), and $U = O(\gamma^{-1})$ is small.  These regions are named in reference to solutions in phase space, where the outer region is the part of the phase space close to the $(U,W)$ plane, while the inner region corresponds to the neighbourhood of the $V$-axis, on which $(U,W)= (0,0)$.

The asymptotic structure of our system is very similar to that which occurs when the Fisher--KPP system is modified by a small cut-off factor in the reaction term \eqref{eq:fkppCutoff}, as considered in \citet{brunet1997shift} and \citet{dumortier2007critical}.  In our case, $\gamma^{-1}$ plays the role of the small parameter $\epsilon$ of that model.  The relationship between $\gamma$ and the wavespeed correction $\delta$ is found by matching the solutions in the intermediate and inner regions.  To leading order this is universal (indeed, the same as \eqref{eq:fkppCutoffSpeed}) and does not depend on the value of $V_\infty$.  However, we will also numerically determine the prefactors in the leading order asymptotic approximations, which are required to observe the effect of $V_\infty$ on the asymptotic wave speed.  

\subsection{Outer and intermediate regions}

\begin{figure}
\centering
\includegraphics[width=\textwidth]{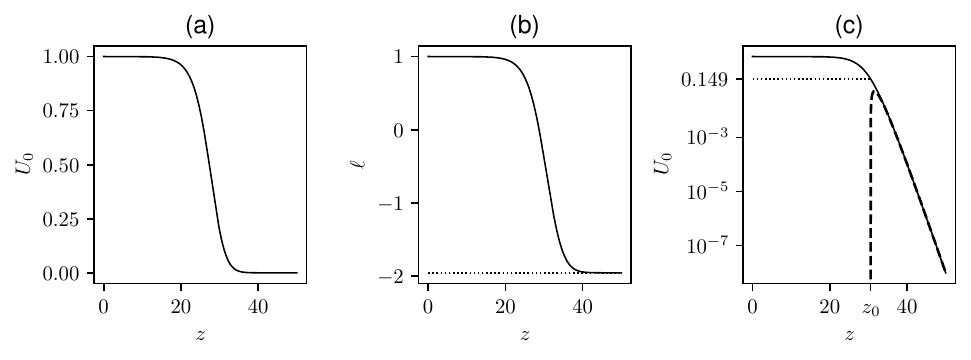}
\caption{(a) Numerically computed solution $U_0(z)$ of the leading-order outer problem \eqref{eq:FisherKppSystem}, which is the Fisher-KPP system.  (b) Plot of the expression $\ell(z)$ \eqref{eq:ell} determined from $U_0$ and $W_0 = U_0'$, which is used to find the prefactor $A_0$ in the $z\to \infty$ limit \eqref{eq:U0AsymptoticBehaviour} of the leading order problem.  (c) Fit of the far field behaviour \eqref{eq:U0AsymptoticBehaviour} with the numerical solution (plotted on a logarithmic scale), showing that when $z$ is equal to the reference value $z_0$, $U_0 \approx 0.149$.  }
\label{fig:U0Asymptotics}
\end{figure}

We first consider the outer region depicted in Fig.~\ref{fig:asymptoticSchematic}.  
In the limit $\gamma \gg 1$, \eqref{eq:Veq} implies that $V$ is exponentially small, so that to leading order, $U$ and $W$ are given by travelling wave solutions to the Fisher--KPP system \eqref{eq:fkpp} with wavespeed $c = 2$:
\begin{equation}
U \sim U_0(z), \qquad W \sim W_0(z), \qquad \delta\to 0,
\label{eq:outerSolution}
\end{equation}
where
\begin{equation}
\frac{\mathrm dU_0}{\mathrm dz} = W_0, \qquad \frac{\mathrm dW_0}{\mathrm dz} = -2W_0 - U_0(1-U_0), \qquad \lim_{z\to-\infty} (U_0, W_0) = (1,0).
\label{eq:FisherKppSystem}
\end{equation}
As it will be required to compute prefactors, we have numerically calculated this solution, using the standard method of solving the initial value problem starting on the unstable eigenvector near $U_0=1, W_0=0$.  A plot of $U_0(z)$ is shown in Fig.~\ref{fig:U0Asymptotics}(a).

As $z\to \infty$, the behaviour of this leading order solution is determined by the repeated eigenvalue with value $-1$ of the Fisher--KPP system \eqref{eq:FisherKppSystem} near the origin.  This repeated eigenvalue implies that to leading order,
\[
U_0 \sim Az\e^{-z}, \qquad z \to \infty
\]
where the prefactor $A$ is arbitrary, as the system \eqref{eq:FisherKppSystem} is translationally invariant in $z$.  However, in order to carry out a matched asymptotic expansion to the order required to observe the effect of different values of $V_\infty$, we require a more precise form of the far field behaviour of $U_0$.  In particular, we write
\begin{equation}
U_0 = A_0(z-z_0)\e^{-(z-z_0)} + o\left(\e^{-z}\right), \qquad z\to \infty,
\label{eq:U0AsymptoticBehaviour}
\end{equation}
where we have explicitly included the translation constant $z_0$.  In this form, the prefactor $A_0$ is not arbitrary, but is a constant that is fixed by stipulating that the order of error be less than exponential.  Since $A_0$ depends on the trajectory of the solution as a whole, it can only be estimated from a numerical solution of \eqref{eq:FisherKppSystem}.

One way to estimate the value of $A_0$ from such a numerical solution of \eqref{eq:FisherKppSystem} is to differentiate \eqref{eq:U0AsymptoticBehaviour} to find $W_0$, and eliminate $A_0$ to find
\[
U_0 + W_0 \sim A_0\e^{-(z-z_0)}, \qquad \frac{U_0}{U_0 + W_0} \sim z-z_0, \qquad z\to \infty.
\]
By taking the logarithm of the first of these, and using the second to eliminate $z-z_0$, we find
\begin{equation}
\log(A_0) = \lim_{z\to\infty} \ell(z), \qquad \ell(z) = \frac{U_0}{U_0+W_0} + \log|U_0+W_0|.
\label{eq:ell}
\end{equation}
In Fig.~\ref{fig:U0Asymptotics}(b) we plot $\ell(z)$ over time, and observe that it tends to a constant $\approx -1.953$ as $z\to \infty$.  Using this we determine that the prefactor is
\begin{equation}
A_0 \approx 0.1419.
\end{equation}

In writing \eqref{eq:U0AsymptoticBehaviour}, the constant of the translation $z_0$ may be thought of as a reference point in the outer region, corresponding to the point where $U_0$ takes a specific (order one) value.  Indeed, this value may also be estimated from the numerical solution by fitting \eqref{eq:U0AsymptoticBehaviour}, with our computed value of $A_0$, to the numerical solution.  We show this fit in Fig.~\ref{fig:U0Asymptotics}(c).  While the numerical value of $z_0$ is arbitrary, as the problem is translationally invariant, we find that it corresponds to the point where $U_0(z_0) \approx 0.149$.  In terms of the asymptotic structure of the problem, what is important is that $z_0$ is a point in the outer region where $U = O(1)$ is a fixed value.

We now consider the intermediate region in which the effect of subcritical wave speed ($\delta > 0$) becomes important while $V$ is still exponentially small.  As the wave speed is below the critical value $2$, the leading order Fisher--KPP approximation is pushed into the regime in which the origin is a stable spiral.  That is, 
\begin{equation}
U \sim U_M, \qquad (2-\delta)\frac{\mathrm dU_M}{\mathrm dz}+ \frac{\mathrm d^2U_M}{\mathrm dz^2} + U_M = 0,
\end{equation}
with asymptotic solution (matching to the outer solution \eqref{eq:U0AsymptoticBehaviour})
\begin{equation}
U_M \sim A_0\delta^{-1/2} \e^{-(z-z_0)} \sin\left(\delta^{1/2}(z-z_0)\right).
\label{eq:UMAsymptoticBehaviour}
\end{equation}
Here the prefactor $A_0\delta^{-1/2}$ results from matching with the outer solution \eqref{eq:outerSolution}, so that expanding $U_M$ in small $\delta$ results in the leading order behaviour \eqref{eq:U0AsymptoticBehaviour} as $z \to \infty$.  This asymptotic behaviour holds when $\delta \to 0$ with $z-z_0 = O(\delta^{-1/2})$, thus the magnitude of $U_M$ is exponentially small despite the algebraically large prefactor.  

The intermediate solution is valid until the sinusoidal term in \eqref{eq:UMAsymptoticBehaviour} approaches zero.  In this limit, the intermediate solution must be matched to the inner solution where $V$ is no longer negligible.  Noting that $z = z_0$ is not in the intermediate region, the first (that is, the smallest) value of $z$ at which this occurs is when $z = z_0 + \delta^{-1/2}\pi$.  The matching between the intermediate and inner problems is where the dependence of $\delta$ on $\gamma$ will be determined, as we describe next.

\subsection{Inner problem and matching}

In the inner region, $V = O(1)$ is required to satisfy the far field condition $V\to V_\infty$ as $z\to\infty$.  From \eqref{eq:Veq}, this is only possible if $U$ (and thus $W$) is small, in particular of order $\gamma^{-1}$.

Formally, let $U = \gamma^{-1}\hat U$, $W = \gamma^{-1}\hat W$.  To leading order (for large $\gamma$), the system \eqref{eqs:odeSystem} is
\begin{align}
\frac{\mathrm d\hat U}{\mathrm dz}  = \hat W, \quad \frac{\mathrm dV}{\mathrm dz} = \frac{1}{2}\hat U V, \quad \frac{\mathrm d\hat W}{\mathrm dz} = \frac{1}{1-V}\left(\frac{1}{2}\hat U V\hat W - 2\hat W - \hat U(1-V)\right).
\label{eq:innerProblem}
\end{align}
Considering solutions for which $V_\infty = V_s$, the relevant boundary condition in this region is that the trajectory approaches the given $V=V_\infty$ along the more negative eigenvector, that is,
\begin{equation}
\begin{bmatrix}
\hat U \\ V \\ \hat W 
\end{bmatrix}
\sim 
\begin{bmatrix}
 0 \\ V_\infty \\ 0\end{bmatrix}
+ \begin{bmatrix} \Lambda \\ V_\infty/2 \\ \Lambda^2 \end{bmatrix} \e^{\Lambda z}, \quad z\to\infty,
\label{eq:innerProblemIC}
\end{equation}
where
\begin{equation}
\Lambda = -\frac{1}{1-V_\infty}\left(1 + \sqrt{1 - (1-V_\infty)^2}\right)
\end{equation}
is the leading order ($c\to 2$) eigenvalue corresponding to $\lambda_3'$ in \eqref{eq:eigsAtInfinity}. 

This inner problem is too intractable to solve in closed form.  However, we only require the far field ($z \to -\infty$) behaviour of these solutions to match to the intermediate region \eqref{eq:UMAsymptoticBehaviour}.  As $z\to -\infty$, $V\to 0$ and the leading order inner problem \eqref{eq:innerProblem} approaches the linearised version of the critical-speed Fisher--KPP equation near the origin, that is
\[
\frac{\mathrm d\hat W}{\mathrm dz} \sim -2\hat W - \hat U, \qquad \frac{\mathrm d\hat U}{\mathrm dz} = \hat W, \qquad z \to -\infty.
\]
Thus, similar to the outer problem, 
\begin{equation}
\hat U \sim -A_I(z-z_I)\e^{-(z-z_I)}, \qquad z\to -\infty.
\label{eq:UIAsymptoticBehaviour}
\end{equation}
As with the analysis of the outer problem, the constant $A_I $ depends in a nontrivial way on the entire solution to \eqref{eq:innerProblem}, and thus on $V_\infty$.  The translation constant $z_I$ may be thought of as a reference location in the inner region, where $\hat U$ takes a specific value (we omit the specific numerical calculation here, as it also depends on $V_\infty$).

\begin{figure}
\centering
\includegraphics[width=\textwidth]{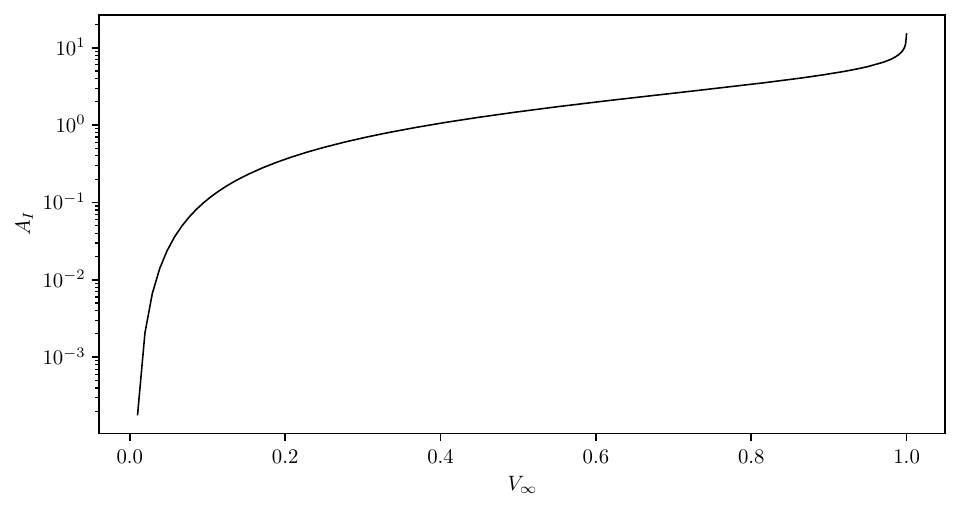}
\caption{Numerically calculated prefactor $A_I$ in the far field behaviour \eqref{eq:UIAsymptoticBehaviour} of the inner problem \eqref{eq:innerProblem}, as dependent on $V_\infty$.  The constant $A_I$ is required in order to perform the matching with the outer problem, and ultimately determine the first correction to the wave speed $c$ as a function of parameter $\gamma$ and $V_\infty$.}
\label{fig:VinfAI}
\end{figure}

With the sign convention in \eqref{eq:UIAsymptoticBehaviour}, $A_I$ will be positive.  By integrating $\eqref{eq:innerProblem}$ numerically in the negative $z$ direction, from an initial condition \eqref{eq:innerProblemIC}, we may numerically compute $A_I$ independently of $z_I$ in a similar way as was done in the outer problem \eqref{eq:ell}:
\[
\log(A_I) = \lim_{z\to -\infty} \left[\frac{\hat U}{\hat U+\hat W} + \log\left\lvert\hat U+\hat W\right\rvert \right].
\]
For the typical values of $V_\infty \in \{0.25, 0.5, 0.75\}$ we have considered thus far, we use this expression to estimate
\begin{equation}
\left.A_I\right\rvert_{V_\infty = 0.25} \approx 0.515, \qquad \left.A_I\right\rvert_{V_\infty = 0.5} \approx 1.485, \qquad \left.A_I\right\rvert_{V_\infty = 0.75} \approx 2.943.
\end{equation}
For completeness, in Fig.~\ref{fig:VinfAI} we plot numerically calculated values of $A_I$ over a range of $V_\infty$.

\begin{figure}
\centering
\includegraphics[width=\textwidth]{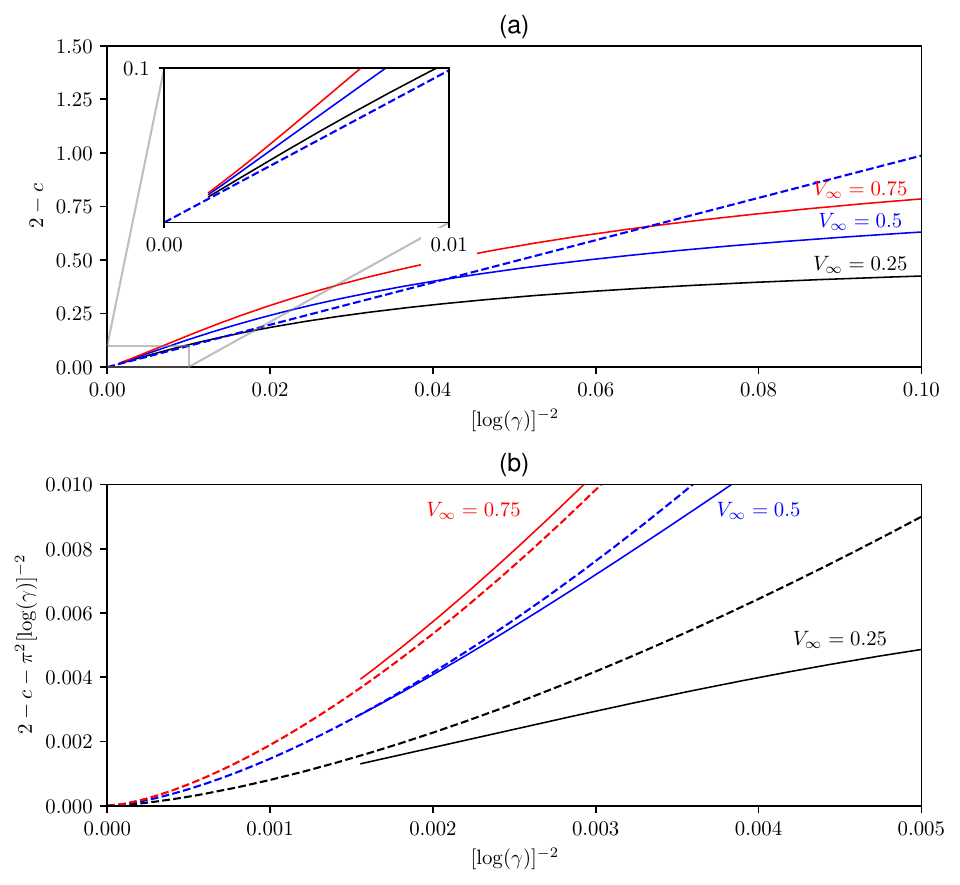}
\caption{Comparison between the numerically determined travelling wave speeds and the asymptotic approximation for large $\gamma$.  (a) Comparison between the numerically determined wavespeed correction $\delta = 2-c$ (solid lines) and the leading order term in \eqref{eq:asymptoticDelta} (dashed line) show convergence to the leading order universal logarithmic behaviour that is independent of $V_\infty$.  (b) The difference between the numerically determined wave speed correction and the leading order approximation (solid lines), compared with the next term in the asymptotic approximation \eqref{eq:asymptoticDelta} (dashed lines), in which the effect of $V_\infty$ plays a role.  These results are calculated up to $\gamma = 10^{11}$, at which point $[\log(\gamma)]^{-2} \approx 0.0015$.}
\label{fig:asymptoticWaveSpeed}
\end{figure}

We now carry out the matching between the intermediate and inner solutions.  In order for the inner solution \eqref{eq:UIAsymptoticBehaviour} to match to the intermediate solution \eqref{eq:UMAsymptoticBehaviour}, the sinusoidal term in the intermediate solution must approach zero, which first happens when its argument is in the neighbourhood of $\pi$.  We define a new spatial variable $\xi$ by $z - z_0 = \delta^{-1/2}\pi + \xi$.  Then expanding the intermediate solution in small $\delta$ with $\xi=O(1)$ we obtain
\[
U_M \sim A_0\delta^{-1/2}\e^{-\delta^{-1/2}\pi - \xi}\sin(\pi + \delta^{1/2}\xi) \sim -A_0\e^{-\delta^{-1/2}\pi} \xi\e^{-\xi}.
\]
To match this expression with \eqref{eq:UIAsymptoticBehaviour} we require $\xi = z - z_I$, so that
\begin{equation}
z_I - z_0 \sim \delta^{-1/2}\pi.
\label{eq:intermediateRegionSize}  
\end{equation}
Equating the prefactors between $U_M$ and the inner solution $\gamma^{-1}\hat U$, we find
\[
-A_0\e^{-\pi\delta^{-1/2}} \sim -A_I\gamma^{-1},
\]
from which we obtain the relationship between the wave speed correction $\delta$ and $\gamma$:
\begin{equation}
\delta \sim \pi^2\left[\log\left(\frac{A_0}{A_I}\gamma\right)\right]^{-2} \sim \frac{\pi^2}{[\log\gamma]^2} + \frac{2\pi^2\log(A_I/A_0)}{[\log\gamma]^3}.
\label{eq:asymptoticDelta}
\end{equation}
The first term on the right hand side of \eqref{eq:asymptoticDelta} is equivalent to the universal correction that arises in the Fisher--KPP equation with cut-off \citep{brunet1997shift,dumortier2007critical}.  It does not depend on the specific value of $V_\infty$, and indeed we expect it would not depend on the specifics of the diffusion or reaction terms being considered, as it is purely a consequence of the fact that the intermediate problem has to match to an inner problem where $U$ is of order $\gamma^{-1}$.  In Fig.~\ref{fig:asymptoticWaveSpeed}a we present numerical evidence that this leading order term is correct, by plotting $(\log\gamma)^{-2}$ against the correction to the velocity $\delta = 2-c$, for the travelling wave branches we have computed for each of $V_\infty \in\{ 0.25, 0.5, 0.75 \}$.  Each of these curves is seen to approach the line with slope $\pi^2$, as predicted by the leading order term in \eqref{eq:asymptoticDelta}.  The convergence to this limit is very slow, which is unsurprising given that the asymptotic series is logarithmic in $\gamma$.

To test the correction term (the second term on the right hand side of \eqref{eq:asymptoticDelta}) we compute the difference between the velocity correction and the leading order term, that is
\[
\delta_1 = 2 - c - \frac{\pi^2}{[\log\gamma]^2},
\]
and, in Fig.~\ref{fig:asymptoticWaveSpeed}b, plot this quantity as a function of $(\log\gamma)^{-2}$, along with the asymptotic prediction.  The asymptotic prediction depends on $V_\infty$ through $A_I$, so is different for each branch.  This comparison confirms that the correction term improves the estimate in the limit $\gamma\to\infty$, but again highlights the very slow convergence of the asymptotic approximation in this limit.

Finally, we note that in the intermediate region in which $U$ behaves as \eqref{eq:UMAsymptoticBehaviour}, both $U$ and $V$ are small.  This region therefore corresponds to the interstitial region that can be observed in Fig.~\ref{fig:bigGammaProfiles}.  Combining \eqref{eq:intermediateRegionSize} with the asymptotic expression for $\delta$ \eqref{eq:asymptoticDelta} predicts that the width of this region to leading order is
\begin{equation}
z_I-z_0 \sim \log(\gamma)
\label{eq:intermediateRegionSizeInGamma},
\end{equation}
hence the slow growth of this region in the large $\gamma$ limit.

\section{Discussion}
\label{sec:discussion}

In this article we have determined the properties of the three dimensional phase space that select the wave speed $c$ of the system \eqref{eqs:PDEs} for compactly supported or exponentially decaying initial conditions.  Although the phase-space properties of the three-dimensional system in travelling wave coordinates \eqref{eqs:odeSystem} are more complicated than the two-dimensional Fisher--KPP system, the general principle of wave speed selection is the same.  For an initial condition of given far field exponential decay, a PDE solution will evolve toward the travelling wave solution which has the same decay rate, if such a nonnegative solution exists, leading to the formula \eqref{eq:speed1}.  Otherwise, the solution will tend to the travelling wave solution which has the largest decay rate.  Depending on the values of $\gamma$ and $v_0$, the solution with fastest decay rate will either correspond to the wave speed $c$ such that $v_0 = V_c$, where $V_c$ is the threshold between monotonic and oscillatory decay at the front \eqref{eq:speed2}, or such that $v_0 = V_s$, where $v_s$ separates exponential decay from above and below, if that point exists.  An important point regarding the determination of $V_s$ is, although it is a property of the phase space, it is not a local property in the neighbourhood of the $V$-axis; rather, it is related to the manner in which the two-dimensional unstable manifold of $(U,V,W) = (1,0,0)$ intersects the $V$-axis.  For this reason, the existence or location of $V_s$ cannot be determined by a local expansion, so there is unlikely to be a general closed-form expression for the selected wave speed.

We have also performed the asymptotic analysis of travelling waves in the large $\gamma$ limit.  In this limit we have focused on solutions selected by $V = V_s$, as particular relevant for compactly supported initial conditions, although this is the appropriate value for $\gamma\to\infty$ for any initial condition with sufficiently large decay such that the wave speed $c < 2$ (that is, the constraint \eqref{eq:stupidSpeed} is not satisfied).  In carrying out the $\gamma\to \infty$ limit it is useful to think of the system \eqref{eqs:odeSystem} as a generalisation of the Fisher--KPP system in which the reaction term is perturbed, leading to a logarithmically small correction to the wave speed from the critical Fisher--KPP wave speed of $2$.  The leading order correction in this case is universal, so is likely to be the correction observed for a large class of systems that look like \eqref{eqs:PDEs} with variations in diffusivity and reaction terms; for example, the system considered by \citet{colson2021travelling}, in which the resident cell population effects the diffusivity but not the reaction, is likely to have the same leading order asymptotic behaviour.

More generally, the universality of the perturbed Fisher--KPP asymptotic structure means that it plays a role in the asymptotic structure of more complex models, such as the \citet{gatenby1996reaction} model \eqref{eqs:gg}.  The logarithmic dependence of gap width on cell degradation rate certainly appears universal.  Specific details such as prefactors can depend on the distinguished limits considered, however.  The leading-order asymptotic approximation for the interstitial gap width ($z_+$) found in \citet{gatenby1996reaction} and \citet{fasano2009slow}, and established rigorously in \cite{davis2022traveling}, is
\[
z_+ \sim \frac{1}{\sqrt{\gamma_2}}\log\alpha, \qquad d \to 0,
\]
where we have used the nondimensional variables in the scaling \eqref{eqs:gg}.  Importantly, the limit in which this holds is vanishing diffusivity ratio $d\to 0$, where $\gamma_2$ is taken to be $O(d)$; the prefactor of $1/\sqrt{\gamma_2}$ then arises from the balance between the reaction and diffusion terms for the acid concentration in \eqref{eqs:gg}.  In contrast, the two-component models such as that considered in our study are derived by assuming $\gamma_2 \gg 1$, so that the invading cell population and acid concentration are effectively coupled.  The logarithmic dependence of interstitial gap separation on degradation rate is maintained, however.

From a biological perspective, the results of the above asymptotic analysis imply that when the interstitial gap is appreciable, there is only very weak connection between the interstitial gap width and the wave speed.  
The invading wave speed is essentially set by the invasion dynamics of the cancer cells in the absence of any resistance by crowding of resident cells (that is, approximately the Fisher--KPP wave speed).  In this case, the wave speed is related to the invading cell parameters only (that is, the invading cell diffusivity and proliferation rates), while the interstitial gap width provides only order-of-magnitude information on the effect of invading cells on the environment, including both the excess acidity, and the effect of that acidity on resident cells and extracellular matrix.  Thus, while the presence of an interstitial gap may indicate more aggressive tumour growth, there is little difference between different width gaps if they are present.

In this work we have not considered the degenerate-diffusion case where the initial resident cell population is $v_0 = 1$ (or $V_\infty = 1$ in calculating travelling wave solutions), at which point the diffusivity of invading cells vanishes.  In phase space such a trajectory would have to either approach $(U,V,W) = (0,1,0)$ as $z\to \infty$ (if solutions remain smooth), or a non-smooth condition where $(U,V) \to (0,1)$ as $z$ tends to some finite value, with $W$ either approaching a constant or becoming unbounded.  For the first such case, it can be shown by substituting an algebraic expansion into \eqref{eqs:odeSystem} that there is a solution that behaves like
\[
U \sim \frac{2c^2}{\gamma} z^{-2}, \qquad V \sim 1 - 2cz^{-1}, \qquad W \sim -\frac{4c^2}{\gamma}z^{-3}, \qquad z\to \infty.
\]
However, such a solution decays algebraically as $z\to\infty$, which is not the behaviour seen in the numerical simulations; instead, one would expect degenerate diffusion to result in a solution that is non-smooth at a moving front \citep{el2021travelling}.  We have not been able to ascertain a consistent expansion of the system \eqref{eqs:odeSystem} with such behaviour, which we leave for future exploration.

\paragraph{Data availability} Data sharing not applicable to this article as no datasets were generated or analysed during the current study.

\subsection*{Declarations}

\paragraph{Conflict of interest} The authors declare that they have no conflicts of interest.

\paragraph{Ethical statement} This study did not involve any experiments with human participants or animals.

\bibliography{refs}

\end{document}